\documentclass[a4paper,12pt,leqno]{article}
\usepackage{latexsym}
\usepackage[all]{xy}

\usepackage{amssymb} 
\usepackage{amsmath} 
\usepackage{theorem}

\def\Z{{\mathbb{Z}}}
\def\R{{\mathbb{R}}}
\def\A{{\mathcal{A}}}

\DeclareMathOperator{\Der}{Der}

\numberwithin{equation}{section}

\newcommand{\owari}{\hfill$\square$}

\theoremstyle{break}
\newtheorem{theorem}{Theorem}[section]

\newtheorem{proposition}[theorem]{Proposition}

\newtheorem{lemma}[theorem]{Lemma}
\newtheorem{define}[theorem]{Definition}

\title{
Primitive filtrations of the modules of
invariant logarithmic forms of
Coxeter arrangements
}
\author{Takuro Abe\thanks{
Department of Mathematics, Kyoto University, 
Kitashirakawa-Oiwake-cho, Sakyo-Ku, 
Kyoto 606-8502, Japan.
email:abetaku@math.kyoto-u.ac.jp.} and 
Hiroaki Terao\thanks{
Department of Mathematics, Hokkaido University, 
Kita-10, Nishi-8, Kita-Ku, 
Sapporo, Hokkaido 060-0810, Japan.
email:terao@math.sci.hokudai.ac.jp.}
}
\date{\today} 

\begin{document}

\maketitle

\begin{abstract}
We define  {\bf primitive derivations}
for Coxeter arrangements which may not be
irreducible.
Using those derivations, we introduce the
{\bf
primitive filtrations}
of the module of invariant logarithmic differential forms
for an arbitrary Coxeter arrangement with an arbitrary
multiplicity.
In particular, when
the Coxeter arrangement is irreducible with a
constant multiplicity,
the primitive filtration
was studied in \cite{AT08}, which
generalizes the Hodge filtration 
introduced by K. Saito (e.g., \cite{S03}).
\end{abstract}

\section{Introduction}
Let $V$ be an $\ell$-dimensional Euclidean space
and
$\A$ be an arrangement of hyperplanes
in $V$. 
We use \cite{OT} as a general reference for arrangements.
For each $H\in\A$, choose a linear form $\alpha_{H} \in V^{*} $
such that $\ker(\alpha_{H} ) = H$.   
Their product $Q := \prod_{H\in\A} \alpha_{H} $
lies in 
the symmetric algebra
$S:=\mbox{Sym}^*(V^*)$.
The quotient field of $S$ is denoted by $F$. 
Let 
$\Omega_S$
and 
$\Omega_F$
denote the $S$-module of regular differential $1$-forms on $V$
and the $F$-vector space of rational differential
$1$-forms on $V$ 
respectively. 
Define 
the $S$-module 
$\Omega(\A,\infty)$
of 
{\bf
logarithmic 
differential $1$-forms}
by 
\begin{multline*}
\Omega(\A,\infty):
=\{\omega \in \Omega_F ~|~Q^{N} \omega 
\mbox{~and~}
(Q/\alpha_{H})^{N} d \alpha_H\wedge \omega 
\mbox{ ~are both regular}\\
\mbox{ for all }H 
\in \A
\mbox{~for~}N \gg 0
\}.
\end{multline*}
In other words, $\Omega(\A, \infty)$ consists of
all logarithmic differential $1$-forms in the sense of 
Ziegler \cite{Z89}.  

Suppose that $\A$ is a {\bf 
Coxeter arrangement}.
Then the corresponding {\bf
Coxeter group} 
$W = W(\A)$ 
naturally acts on $V$, $V^{*}$,  $S$ and $\Omega(\A, \infty)$.
%
Note that we do not assume that $\A$ is irreducible.
When $\A$ is irreducible, the primitive derivations 
play the central role to define the
Hodge filtration introduced
by K. Saito.  (See \cite{S03} for example.)
In this paper we develop a theory of primitive
derivations and the Hodge filtration in the case of 
non-irreducible Coxeter arrangements.  More precisely,
in Section \ref{Proof}, we introduce primitive derivations 
even when $\A$ is not irreducible. 
Fix a primitive derivation $D$.
Let $R := S^{W} $ be the
$W$-invariant subring of $S$ and
$$T := \{
f\in R ~|~ D(f) = 0
\}.$$
Consider the $T$-linear connection 
(covariant derivative) 
\[
\nabla_{D} : 
\Omega_{F} 
\to
\Omega_{F} 
\]
characterized by 
(1)
$\nabla_{D} (f\omega)
=
D(f) \,\omega + f\, (\nabla_{D}\, \omega)$ 
for $f\in F$ and $\omega\in\Omega_{F} $
and
(2)
$\nabla_{D} (d\alpha) = 0$ 
for all $\alpha\in V^{*}$.
Our first main result is

\begin{theorem}
Let $\Omega(\A, \infty)^{W} $ be the $W$-invariant part of 
$\Omega(\A, \infty)$.
Then the 
$\nabla_{D} 
$ 
induces
a $T$-linear automorphism
\[
\nabla_{D} 
:
\Omega(\A, \infty)^{W} \
\widetilde{\longrightarrow}
\
\Omega(\A, \infty)^{W}. 
\]
\label{thm:automorphism}
\end{theorem}

Note that the inverse map
$\nabla_{D}^{-1}$
and $\nabla_{D}^{k} \ (k\in \Z) $ 
are also $T$-automorphisms. 
Under the assumption that $\A$ is irreducible,
 Theorem \ref{thm:automorphism} 
was proved in \cite[Theorem 1.2 (1)]{AT08}. 

\begin{define}
\label{def:primitivefiltration}
Let $I^{*} : \Omega_{F}  \times \Omega_{F} \rightarrow F$ be 
the $F$-bilinear map induced from the inner product $I$ 
of the Euclidean space $V$. 
Let 
$\bf m : \A \rightarrow \Z$ be an arbitrary multiplicity.
Define
\begin{multline*} 
\Omega(\A, {\bf m})
:=
\{
\omega\in \Omega(\A, \infty)
~|~
(Q/\alpha_{H})^{N} \alpha_{H}^{{\bf m}(H)} I^{*}(d\alpha_{H}, \omega ) \in S \\
\mbox{~for all ~} H \in \A 
\mbox{~for~} N \gg 0 
\}
\end{multline*} 
and
$$
\Omega(\A, {\bf m})^{W}
:=  
\Omega(\A, {\bf m}) \cap 
\Omega(\A, \infty)^{W}.
$$
The {\bf
primitive filtration
of
$\Omega(\A, \infty)^{W}$
 induced from ${\bf m}$} 
is given by
\[
P^{({\bf m})}_{k} 
:=
\nabla_{D}^{k} 
\Omega(\A, {\bf m})^{W}
  \ \ (k\in\Z).   
\]
\end{define} 
Note that 
\begin{multline*} 
\Omega(\A, {\bf m})
=
\{
\omega\in \Omega_{F} 
~|~
(\prod_{H\in\A} \alpha_{H}^{{\bf m}(H) })\omega 
\mbox{~and~}
(\prod_{H\neq H_{0} } \alpha_{H}^{{\bf m}(H) } )(
d\alpha_{H_{0} }  \wedge\omega) \\
\mbox{~are both regular~for all~} H_{0}  \in \A 
\}
\end{multline*} 
if ${\bf m}(H) \ge 0$ for all $H\in\A$.
In this case, $\Omega(\A, {\bf m})$ was
introduced by Ziegler \cite{Z89}.

Our second main result is
an explicit description of the primitive filtration:

\begin{theorem}
The primitive filtration is an increasing filtration
\[
\dots \subset P^{({\bf m})}_{-1}  \subset P^{({\bf m})}_{0}
  \subset P^{({\bf m})}_{1} \subset \dots
\]
such that 
\[
P^{({\bf m})}_{k} = P^{({\bf m} + 2k)}_{0}  
=
\Omega(\A, {\bf m} + 2k)^{W} 
\]
where $({\bf m} + 2k)(H) = {\bf m}(H) + 
2k \ (k\in\Z, H\in\A).$ 
\label{thm:primitivefiltration}
\end{theorem}
When $\A$ is irreducible and
${\bf m}$ is equal to
the constant function 
${\bf 1} $ with ${\bf 1}(H) = 1
 \ (H\in\A)$, the primitive filtration 
coincides with the filtration introduced in  
\cite{AT08}.
Its dual version in Theorem \ref{thm:primitivefiltrationD} 
generalizes the Hodge filtration
introduced by K. Saito (e.g., \cite{S03}).

We construct bases for the primitive filtration induced from 
$\bf 1$ in Theorem \ref{thm:bases}.
The bases are used when we
 prove Theorems \ref{thm:automorphism}  
and \ref{thm:primitivefiltration} in Section \ref{sec:3}.

In Section \ref{sec:4}, we translate our main results
Theorems \ref{thm:automorphism},
 \ref{thm:primitivefiltration} 
and \ref{thm:bases} into the dual language in terms of
the logarithmic derivations. 

\section{Primitive derivations}
\label{Proof} 

We first state a multiple version of
Saito's criterion due to Abe 
\cite{A08}.

\begin{proposition}
\label{saitoziegler} 
Let $\A$ be a central arrangement in $V$ with an
arbitrary multiplicity 
${\bf m} : \A \rightarrow \Z$.
Let $x_{1} , x_{2} , \dots , x_{\ell} $ be a basis for
$V^{*} $.
Define $$Q^{\bf m} := 
\prod_{H\in \A}  \alpha_{H}^{{\bf m}(H)}\in F.  $$ 
Let
$\omega_{1} , \omega_{2} , \dots, \omega_{\ell} 
\in 
\Omega(\A, {\bf m})$.   
Then

(1) $Q^{\bf m} (\omega_{1} \wedge  \dots\wedge \omega_{\ell})$ is regular,

(2) 
$\omega_{1} , \omega_{2} , \dots, \omega_{\ell}$ form an $S$-basis for 
$\Omega(\A, {\bf m})$ if and only if   
$$Q^{\bf m} (\omega_{1} \wedge  \dots\wedge \omega_{\ell}) 
\in
\R^{\times} 
(dx_{1} \wedge dx_{2} \wedge \dots\wedge dx_{\ell} ).
$$
\end{proposition}

\noindent
{\em Remark.}
When ${\bf m} = {\bf 1}$, this is due to K. Saito \cite{S80}.
When ${\bf m}: \A \rightarrow \Z_{\ge 0} $, this is due to Ziegler
\cite{Z89}.   

\medskip

The original proof in 
\cite[Theorem 1.4]{A08}
is written in a slightly 
different language from this paper, so
we include our proof here.

\medskip

\noindent
{\em Proof of Proposition \ref{saitoziegler}.} 
Pick $H\in\A$ arbitrarily and fix it.
Let $m = {\bf m}(H). $
Choose an orthonormal basis 
$x_{1} , x_{2} , \dots, x_{\ell} $ such that
$H = \{x_{1} = 0\}$.

(1) It is enough to show that 
$
x_{1}^{m}
(\omega_{1} \wedge  \dots\wedge \omega_{\ell}) 
$   has no pole along $H$.
Write 
\[
\omega_{j} = \sum_{i=1}^{\ell} f_{ij} dx_{i}
\,\,\,\,
(j=1, \dots , \ell).  
\]
Since $\omega_{j} \in \Omega(\A, {\bf m} )$,
$
x_{1}^{m} f_{1j} 
=
x_{1}^{m}
I^{*} (dx_{1}, \omega_{j})
 $   has no pole along $H$ for 
$
j=1, \dots , \ell
$ by Definition \ref{def:primitivefiltration}.
Moreover,
$$
\sum_{i\ge 2} f_{ij} dx_{1} \wedge dx_{i} =
dx_{1} \wedge \omega_{j} 
$$  
has no pole along $H$ because $\omega_{j} \in \Omega(\A, \infty)$
for  
$j=1, \dots , \ell$.
This implies that $f_{ij} $ has no pole along
$H$  if $i\ge 2$.
Therefore
\[
x_{1}^{m}
(\omega_{1} \wedge  \dots\wedge \omega_{\ell}) 
=
\begin{vmatrix}
x_{1}^{m} f_{11} & x_{1}^{m} f_{12} & \dots & x_{1}^{m} f_{1\ell}\\
f_{21} & f_{22} & \dots & f_{2\ell}\\
\vdots &  \vdots      &  \ddots     &   \vdots \\
f_{\ell 1} & f_{\ell 2} & \dots & f_{\ell\ell}
\end{vmatrix}
(dx_{1} \wedge dx_{2} \wedge \dots\wedge dx_{\ell} )
\]
 has no pole along $H$.

(2)
Suppose that 
 $\omega_{1} , \omega_{2} , \dots, \omega_{\ell}$ form an $S$-basis for 
$\Omega(\A, {\bf m})$.
By (1) we may write
\[
Q^{\bf m} 
(\omega_{1} \wedge \dots \wedge \omega_{\ell} )
=
f 
(dx_{1} \wedge dx_{2} \wedge \dots\wedge dx_{\ell} )
\]
with $f\in S$.
In order to prove that $f$ is a nonzero constant, it is enough to show that
$f$ is not divisible by $x_{1} $.
Define a multiplicity ${\bf m'} : \A \rightarrow \Z_{\ge 0} $ 
by
\[
{\bf m'}(K) :=
\begin{cases}
|{\bf m}(K) |   & \mbox{if~ } K \neq H\\
0   & \mbox{if~ } K = H.
\end{cases}  
\]
Then it is not hard to see that
\[
\eta_{1} := Q^{{\bf m'} } (dx_{1} /x_{1}^{m}), \ \
\eta_{2} := Q^{{\bf m'} } dx_{2}, \ \
\dots, \ \
\eta_{\ell} := Q^{{\bf m'} } dx_{\ell}
\]
lie in $\Omega(\A, {\bf m} )$. 
  Thus
\begin{align*} 
&~~~~(Q^{\bf m}/x_{1}^{m}) (Q^{{\bf m'} } )^{\ell}
(dx_{1} \wedge dx_{2} \wedge \dots\wedge dx_{\ell} )\\
&=
Q^{\bf m} (\eta_{1} \wedge \eta_{2} \wedge \dots \wedge \eta_{\ell})
\in 
S
(\omega_{1} \wedge \dots \wedge \omega_{\ell} )
=    
S f 
(dx_{1} \wedge dx_{2} \wedge \dots\wedge dx_{\ell} ).
\end{align*} 
This implies that $g := (Q^{\bf m}/x_{1}^{m}) (Q^{{\bf m'} } )^{\ell}$
is divisible by $f$.
Since $g$ is not divisible by $x_{1}$, 
neither is $f$.

Suppose that 
$
Q^{\bf m} 
(\omega_{1} \wedge \dots \wedge \omega_{\ell}) 
=
dx_{1} \wedge\dots\wedge dx_{\ell}. 
$ 
In order to prove that 
$
\omega_{1}, \dots , \omega_{\ell} 
$
form a basis it is enough to show that 
$\omega_{1}, \dots , \omega_{\ell} 
$
span $\Omega(\A, \bf m)$ over $S$.
Fix $\omega\in\Omega(\A, \bf m)$.
By (1) we may write
\[
Q^{\bf m} 
(\omega_{1} \wedge \dots \wedge \omega_{i-1} 
\wedge \omega \wedge \omega_{i+1} \dots \wedge \omega_{\ell}) 
=
f_{i} 
(dx_{1} \wedge\dots\wedge dx_{\ell}), 
\]
with $f_{i} \in S$ for $i=1,\dots,\ell$.
Define 
$
\eta := 
\omega - 
\sum_{i=1}^{\ell} f_{i} \omega_{i}.  $  
Then we obtain
\[
Q^{\bf m} 
(\omega_{1} \wedge \dots \wedge \omega_{i-1} 
\wedge \eta \wedge \omega_{i+1} \dots \wedge \omega_{\ell}) 
=
0
\,\,\,
(i = 1,\dots, \ell).
\]
Since $
\omega_{1}, \dots , \omega_{\ell} 
$ 
span the cotangent space of $V$ at each point outside
the hyperplanes, we have $\eta = 0$ and thus
$
\omega = 
\sum_{i=1}^{\ell} f_{i} \omega_{i}.
$  
$\square$

\medskip

Next 
let $\A$ be an irreducible Coxeter arrangement.
Then we may put 
$$
R=S^{W}=\R[P_{1},\ldots,P_{\ell}]
$$ 
with $$
\deg P_{1} < \deg P_{2} \le \cdots \le \deg P_{\ell-1}< 
\deg P_{\ell}
$$
by
\cite{C55}. 
The derivation
\[
D := \frac{\partial}{\partial P_{\ell}} 
\]
is called
a {primitive derivation}
which was extensively studied by K. Saito.
Although $D$ depends upon the choice of 
$P_{\ell}$, its ambiguity  is only
up to a constant multiple.
Recall the $T$-linear connection 
\[
\nabla_{D} : 
\Omega_{F} 
\to
\Omega_{F}. 
\]
Then 
the
$
\nabla_{D} 
$ 
induces 
a $T$-linear automorphism
\[
\nabla_{D} 
:
\Omega(\A, \infty)^{W} \
\widetilde{\longrightarrow}
\
\Omega(\A, \infty)^{W}
\]
by 
\cite[Theorem 1.2 (1)]{AT08}. 
Recall

\begin{proposition} 
\cite[Theorems 1.1 and 2.12]{AT08} 
Suppose that 
$\A$ 
is an irreducible Coxeter arrangement. 
For any $k\in\Z$ and $1\le j\le \ell$ , define
$$
\theta_{j}^{(k)}
:=
\nabla_{D}^{k} \left(dP_j\right)
,
\ \
\Theta^{(k)} := \{\theta_{j}^{(k)}\}_{1\le j\le \ell}
, \ \ \mbox{~and~} \ \
\Theta := \bigcup_{k\in\Z} \Theta^{(k)}.
$$
Then
\begin{itemize}
\item[(1)]
the $S$-module
$
\Omega(\A, 2k-1)
$
is free with a basis
$
\Theta^{(k)}, 
$ 
\item[(2)]
the $R$-module
$
\Omega(\A, 2k-1)^{W}
$
is free with a basis
$
\Theta^{(k)}, 
$ 
\item[(3)]
the $T$-module
$\Omega(\A,2k-1)^{W} $
 is free with a basis
$
\bigcup_{p\le k}
\Theta^{(p)},  
$ and
\item[(4)]
the $T$-module
$\Omega(\A, \infty)^{W} $
 is free with a basis
$
\Theta.  
$ 
\end{itemize}
\label{prop:irreducible}
\end{proposition}

\begin{proposition} 
\cite[Lemma 2.3 and Proposition 2.6 (4)]{AT08} 
\label{prop:recursivematrices}
Let $G:=\left[I^{*} (dP_{i}, dP_{j})\right]_{1\le i, j\le\ell}$.
For each $k\in\Z$, there exists
an $\ell\times\ell$-matrix 
$G_{k}$ with entries in $R$ 
such that
\[
\left[
\theta_{1}^{(k)},
\dots,
\theta_{\ell}^{(k)}
\right]
=
\left[
\theta_{1}^{(k+1)},
\dots,
\theta_{\ell}^{(k+1)}
\right]
G_{k},
\]
where 
$G_{k}$ can be expressed as
$
G_{k} 
=
B_{k} G B'_{k}
$
with 
$B_{k}, B'_{k} \in GL_{\ell}(T).$ 
\end{proposition}

From now on assume that $\A$ is 
an arbitrary Coxeter arrangement which may not
be irreducible.
 Then
 one has the following decompositions:
 \begin{align*}
 V &= V[1] \oplus \cdots \oplus V[t], \ \
 \A = \A[1] \times \cdots \times \A[t],\\
 W&=W[1] \times \cdots \times W[t], \ \
 S\simeq 
S[1] \otimes_\R \cdots \otimes_\R S[t],
 \end{align*}
 where each
 $\A[i]$ is an irreducible Coxeter arrangement
 in $V[i]$, 
 $W[i] := W(\A[i])$,
 and
$$
S[i] := S(V[i]^{*})
=
\R[x_{1}[i], \dots , x_{\ell[i]}[i]]
$$
for $i = 1\,\dots, t$.
We naturally regard 
$\A[i]$ as a subarrangement of $\A$, 
$S[i]$ as a subring of $S$, 
and $W[i]$ as a subgroup of $W$. 

Let $1\le i \le t$. 
Let $R[i]$ denote the $W[i]$-invariant subring of
$S[i]$.  
 Let 
 $\ell[i] = \dim V[i]$.
Then we may put
\[
R[i] = \R[
P_{1}[i], \dots , P_{\ell[i]}[i]
]
\]
with
\[
\deg P_{1}[i] < \deg P_{2}[i] \le  \dots  < \deg P_{\ell[i]}[i].
\]
Then
\[
R = S^{W} =
\R[\{ P_{j}[i]   \}_{1\le i\le t, 1\le j\le \ell[i]} ]
\simeq
R[1] \otimes_\R \cdots \otimes_\R R[t].
\]
Thus we may
naturally regard $R[i]$ as a subring of $R$. 
%
%
Let $D[i] : R[i] \rightarrow R[i]$ 
denote a primitive derivation corresponding to
the irreducible Coxeter arrangement $\A[i]$.
  We may
 naturally extend the derivation $D[i]$
 to a derivation $\hat D[i] : R \rightarrow R$
 by $\hat D[i] (f) = 0$ for any
 $f\in R[j] \ (i\neq j)$.

 \begin{define} 
Let $\A$ be a Coxeter arrangement which may not be irreducible.
Then the 
 derivation 
 $$
 D:=
 \sum_{i=1}^{t} \hat D[i]
 : R\rightarrow R
 $$
 is called a {\bf primitive derivation}
 of $W$.
Let $T := \ker (D : R \rightarrow R)$. 
\label{def:primitivederivation}  
 \end{define}

\noindent
{\em Remark.}
The
 primitive derivations defined in
Definition \ref{def:primitivederivation} 
are
 not necessarily homogeneous or 
 unique up to a constant multiple
 unlike the irreducible case.
However, those derivations play a similar role
to irreducible primitive derivations
as we show in this note.

\bigskip

We often write $P[i]$ instead of
$P_{\ell[i]}[i]$ for simplicity.
Then we have

\begin{lemma}
For
$i = 1, \dots, t$,
$R
=
T[P[i]].
$
\label{lem:T}
\end{lemma}

\noindent
{\em Proof.}
It is obvious that 
$\{ P_{j}[i]   \}_{1\le i\le t, 1\le j\le \ell[i]-1}
\subset
T.
$ 
Note $P[j]-P[i]\in T$ 
because
$$D(P[j]-P[i]) = D(P[j]) - D(P[i]) = 1-1=0.$$
Thus
$P[j] =  
(P[j]-P[i])+P[i]
\in T[P[i]].$
\owari

\medskip

\begin{theorem}
For any $k\in\Z$,
$1\le i\le t$  and $1\le j\le \ell[i]$ , define
$$
\theta_{j}^{(k)}[i]:=
\nabla_{D[i]}^{k} \left(dP_j[i]\right).
$$
Let
$$
\Theta^{(k)}[i] := \{\theta_{j}^{(k)}[i]\}_{1\le j\le \ell[i]}, 
\
\Theta^{(k)} := \bigcup_{i=1}^{t}  \Theta^{(k)}[i],
\ \mbox{~and~} \
\Theta := \bigcup_{k\in\Z} \Theta^{(k)}.
$$
Then
\begin{itemize}
\item[(1)]
the $S$-module
$
\Omega(\A, 2k-1)
$
is free with a basis
$
\Theta^{(k)}, 
$ 
\item[(2)]
the $R$-module
$
\Omega(\A, 2k-1)^{W}
$
is free with a basis
$
\Theta^{(k)}, 
$ 
\item[(3)]
the $T$-module
$\Omega(\A,2k-1)^{W} $
 is free with a basis
$
\bigcup_{p\le k}
\Theta^{(p)},  
$ and
\item[(4)]
the $T$-module
$\Omega(\A, \infty)^{W} $
 is free with a basis
$
\Theta.  
$ 
\end{itemize}

\label{thm:bases}
\end{theorem}

\noindent
\textit{Proof.}  
(1) Let $t=2 $ for simplicity.
By Proposition \ref{prop:irreducible} (1), 
$\Theta^{k}[i]$ is an $S[i]$-basis for 
$\Omega(\A[i], 2k-1)$ for each $i$. 
Thus, by Proposition \ref{saitoziegler},
  we have
\begin{align*} 
&~~~~Q^{2k-1} 
\left(
\bigwedge_{j=1}^{\ell[1]} \theta^{(k)}_{j} [1]   
\right)
\bigwedge
\left(
\bigwedge_{j=1}^{\ell[2]} \theta^{(k)}_{j} [2]   
\right)\\
&=
\left(
Q[1]^{2k-1} 
\bigwedge_{j=1}^{\ell[1]} \theta^{(k)}_{j} [1]   
\right)
\bigwedge
\left(
Q[2]^{2k-1} 
\bigwedge_{j=1}^{\ell[2]} \theta^{(k)}_{j} [2]   
\right)\\
&\in
\R^{\times} 
(
dx_{1}[1] \wedge \dots\wedge dx_{\ell[1]}[1]
\wedge
dx_{1}[2] \wedge \dots\wedge dx_{\ell[2]}[2]
),
\end{align*} 
where $Q[i] = \prod_{H\in \A[i]} \alpha_{H} \ \ (i = 1, 2)$. 
This implies (1) because of Proposition \ref{saitoziegler} again.

(2) 
Note that each $\theta_{j}^{(k)}[i]$ is $W$-invariant by definition. 
Thus 
$$
\Theta^{(k)}  \subset \Omega(\A,2k-1)^W.
$$
Since $\Theta^{(k)} $ is linearly independent over $S$
by (1),
so is over $R$.  An arbitrary element of 
$
\Omega(\A,2k-1)^W 
$
can be expressed as a linear combination of
$\Theta^{(k)}$ with coefficients in $S.$ 
Then it is obvious that each of the coefficients lies in $R$.
This shows that    $\Theta^{(k)}$ spans
$
\Omega(\A,2k-1)^W 
$
over $R$. 

\medskip

(3)
Let 
$$
\mathcal T[i] := 
\bigcup_{p\le k}
\Theta^{(p)}[i]
\mbox{
and
}
\mathcal T :=
\bigcup_{i=1}^{t} 
\mathcal T[i]
= 
\bigcup_{p\le k}
\Theta^{(p)}.
$$

\noindent
{\em Step 1. $\mathcal T$ spans 
$\Omega(\A,2p-1)^W 
$ over $T$.}
 
 Since 
\[
\Theta^{(p)}
\subset
\Omega(\A,2p-1)^W 
\subseteq
\Omega(\A,2k-1)^W 
\]
for $p\le k$,
we have $\mathcal T\subset
\Omega(\A,2k-1)^W 
$.
Let 
$\left<\mathcal T\right>_{T} $ 
be 
the submodule of 
$\Omega(\A,2k-1)^W 
$ generated by $\mathcal T$ over $T$. 
%
Let 
\[
T[i] := \ker (D[i] : R[i] \rightarrow R[i])
\]
for each $i$.  Then $T[i]\subseteq T$.  
By Proposition \ref{prop:irreducible} (3)
we know that
$\left<{\mathcal T}[i]\right>_{T[i]} $ 
is closed under the multiplication of 
$R[i]$ for each $i$. In particular,
$$P[i]\cdot
{\mathcal T}[i]
\subset
\left<{\mathcal T}[i]\right>_{T[i]}
\subseteq
\left<{\mathcal T}[i]\right>_{T}
$$
because $P[i] = P_{\ell[i]}[i] \in R[i]$. 
Therefore
$\left<{\mathcal T}[i]\right>_{T} $ 
is closed under the multiplication of 
$R$
because
$R = T[P[i]]$ by Lemma \ref{lem:T}.
Thus
we obtain
$
\left<
{\mathcal T}[i]
\right>_{R} 
 =
\left<{\mathcal T}[i]\right>_{T}
$
for each $i$. 
Therefore
$
\left<{\mathcal T}\right>_{R}
=
\left<{\mathcal T}\right>_{T}.
$
By (2) we have
\[
\Omega(\A, 2k-1)^{W} 
=
\left<
\Theta^{(k)} 
\right>_{R} 
\subseteq
\left<
{\mathcal T}
\right>_{R} 
=
\left<
{\mathcal T}
\right>_{T} 
\subseteq 
\Omega(\A, 2k-1)^{W}. 
\]
Therefore 
$\left<
{\mathcal T}
\right>_{T} 
=
\Omega(\A, 2k-1)^{W}$: 
${\mathcal T}
$
spans
$
\Omega(\A, 2k-1)^{W}$ over $T$. 

\medskip

\noindent
{\em Step 2. $\mathcal T$ is linearly independent over $T$.}

It is enough to show that 
$\mathcal T[i]$ is linearly independent over $T$ for each $i$.
Let $1\le i\le t$. 
%
%
Assume 
\[
\sum_{k\in\Z} 
\left[
\theta_{1}^{(k)}[i], \dots, \theta_{\ell[i]}^{(k)}[i]
\right]
{\bf g}_{k} ={\bf 0}
\]
with 
${\mathbf g}_{k}  =
\left[
g_{k, 1},\ldots
g_{k, \ell[i]}
\right]^T \in T^{\ell[i]},\ k \in \Z$
 such that 
there exist
integers $p$ and $q$ such that $p \le q$, 
${\mathbf g}_{p} \neq 0,{\mathbf g}_{q}  \neq 0$ and 
${\mathbf g}_{k} =0$ for all $k < p$ and $k > q$. 
Then, by Proposition \ref{prop:recursivematrices}  
$$
{\bf 0}
=\sum_{k=p}^q 
\left[
\theta_{1}^{(k)}[i], \dots, \theta_{\ell[i]}^{(k)}[i]
\right] 
\,{\mathbf g}_{k} 
=
\left[
\theta_{1}^{(q)}[i], \dots, \theta_{\ell[i]}^{(q)}[i]
\right] 
\sum_{k=p}^{q}  H_{k} {\bf g}_{k},  
$$
where
\[
H_{q} :=
I_{\ell[i]},
\ \ \
 H_{k} := G_{k} G_{k+1} \dots G_{q-1}
\ \ \ (p\le k < q). 
\]
This
implies that 
$$
{\bf 0}
=
\sum_{k=p}^q H_{k}{\mathbf g}_{k}.
$$
Note that
$H_{k}$ can be expressed as a product of 
$(q-k)$
copies of 
$$G[i]
:=
\left[
I^{*}(dP_{a}[i], dP_{b}[i]) 
\right]_{1\le a, b \le \ell[i]} 
$$ and 
matrices belonging to $\mbox{GL}_{\ell[i]}(T[i])$. 
It is well-known that 
$
D[G[i]]  
=
D[i][G[i]]  
\in GL_{\ell[i]}(T[i]) $ 
\cite[Proposition 2.1]{AT08}.
Thus
$D^{q-p}[H_{k} ] = 0
\ (k > p)$.
Applying $D^{q-p}$ to the above, we
obtain
$$
D^{q-p}[H_{p}]{\mathbf g}_{p} =0.
$$
Note 
Since the matrix
$D^{q-p}[H_{p}]$,
which is 
a product of 
$(q-p)$ copies of $D[G[i]]$ and 
matrices in $\mbox{GL}_{\ell[i]}(T[i])$,
is
nondegenerate, 
we get
${\mathbf g}_{p} =0$, which is a contradiction. 
This implies that ${\mathcal T}$ is linearly independent
over $T$. 

\medskip

(4) It follows from (3) and the fact that 
$$
\Omega(\A, \infty)^{W} 
=
\bigcup_{k\in\Z} \Omega(\A, 2k-1)^{W}. 
$$ 
\owari
\medskip

%
%

\section{Proof of main theorems.} 
\label{sec:3}

\noindent
{\em Proof of Theorem \ref{thm:automorphism}. }
Since 
\[
\nabla_{D} \theta^{(k)}_{j} [i] 
=
\nabla_{D} \nabla_{D[i]}^{k}  (dP_{j}[i]) 
=
\nabla_{D[i]}^{k+1}  (dP_{j}[i]) 
=
\theta^{(k+1)}_{j} [i], 
\]
the connection $\nabla_{D} $ induces 
a bijection of $\Theta$ to itself.  
Thus $\nabla_{D} $ induces a $T$-automorphism
of $\Omega(\A, \infty)^{W} $ because of Theorem 
\ref{thm:bases} (4). 
\owari

\medskip

For $f\in F$ with $f\neq 0$ 
and $\alpha\in V^{*} \setminus \{{\bf 0} \}$ define
\[
\mbox{ord}_{\alpha} (f)
:=
\min\{
k\in\Z ~|~ \alpha^{k} f \in S_{(\alpha)} 
\},
\]
where $S_{(\alpha)} $ is the localization of $S$ at the prime ideal
$(\alpha) = \alpha S$.
In other words $\mbox{ord}_{\alpha}(f)$  is the order of poles of $f$
  along the hyperplane
$\ker(\alpha)$.

\begin{lemma}
\label{lemma:ord}
Assume that $\A$ is a Coxeter arrangement which may not be irreducible.
Let $D$ be a primitive derivation of $\A$.
Choose $\alpha\in V^{*} $ such that $\ker(\alpha)\in\A$.  
Then

(1) 
$\mbox{\rm ord}_{\alpha} D(\alpha) = 1$.

(2) For $f\in F\setminus\{0\}$ 
with $\mbox{\rm ord}_{\alpha}(f)\neq 0$,
$
\mbox{\rm ord}_{\alpha} (D(f)) 
=
\mbox{\rm ord}_{\alpha} (f) 
+2.
$
\end{lemma}

{\em Proof.}
(1)
Assume that $$\A = \A[1] \times \dots \times \A[t]$$
such that each $\A[i]$ is irreducible. 
Suppose $\ker(\alpha)\in \A[k]$.  Then
$D[i](\alpha) = 0$ if $i\neq k.$  
This implies that we may assume that $\A$ is irreducible
from the beginning.  Choose an orthonormal basis
$\alpha=x_{1}, x_{2}, \dots , x_{\ell}   $
and let $h_{j} := D(x_{j} )$ for $1\le j\le \ell.$  
It is well-known 
(e.g., \cite[pp. 249-250]{ST98} )
that $h_{j} \ (j > 1)$ has no poles along $x_{1} = 0.$
On the other hand,
it is also known 
(e.g., \cite[Corollary 3.32]{ST98} )
that $$\det \left[
\partial h_{j}/\partial x_{i}  \right] = c\, Q^{-2} $$ 
for some nonzero constant $c$.
Thus $h_{1} $ should have poles along $x_{1} =0.$   
Since $Q h_{1} = (QD)(x_{1} )$ is regular,
we have
$
\mbox{\rm ord}_{\alpha} D(\alpha) 
=
\mbox{\rm ord}_{\alpha} h_{1}  
= 1$.

(2)
Suppose that
$k := \mbox{\rm ord}_{\alpha} (f) \neq 0$.  Put
$f = g/\alpha^{k} $.  Then 
$g\in S_{(\alpha)} $
 and
$g\not\in \alpha S_{(\alpha)} $.
Compute
\[
D(f) = D(g/\alpha^{k})
=
D(g)/\alpha^{k} - 
k D(\alpha) g/\alpha^{k+1}. 
\]
From (1)
we have
$\mbox{\rm ord}_{\alpha}(D(\alpha)) 
= 1. $
Since
\[
\mbox{\rm ord}_{\alpha}  (D(\alpha) g/\alpha^{k+1}) = k + 2,
\ \ \
\mbox{\rm ord}_{\alpha} (D(g)/\alpha^{k} ) \le k+1,
\]
we obtain
$
\mbox{\rm ord}_{\alpha} (D(f)) = k+2.
$ 
\owari

\bigskip

\noindent
{\em Proof of Theorem \ref{thm:primitivefiltration}. }
It is enough to prove
$ \nabla_{D} \Omega(\A, {\bf m})^{W}
=
 \Omega(\A, {\bf m}+2)^{W}$. 
Let $\omega\in\Omega(\A, \infty)^{W}$
and $\alpha\in V^{*} $ with $\ker(\alpha) \in \A$.

We first verify:
\begin{equation} 
\mbox{\rm ord}_{\alpha} I^{*} (\omega, d\alpha) 
\neq 0.
\label{eqn:2.1} 
\end{equation} 
Let $s_{\alpha} $ be the orthogonal reflection through 
the hyperplane $\ker(\alpha)$.
Since $\omega$ is $W$-invariant, we have 
$
s_{\alpha} (
I^{*}(\omega, d\alpha)
)
=-
I^{*}(\omega, d\alpha).
 $  
Suppose that
$\mbox{\rm ord}_{\alpha} I^{*} (\omega, d\alpha)=0$.
Then, for a sufficiently  large integer $N$,
$$g := 
(Q/\alpha)^{N} I^{*} (\omega, d\alpha) \in S\setminus \alpha S.$$  
On the other hand, we obtain
\[
s_{\alpha} (g)
=
\left(
s_{\alpha}
(Q/\alpha)^{N}\right) s_{\alpha}  \left(I^{*} (\omega, d\alpha)\right)
=
-(Q/\alpha)^{N} I^{*} (\omega, d\alpha)
=-g.
\]
This shows that $g$ is an antiinvariant with respect to the 
reflection group
$\{{\bf 1}, s_{\alpha} \}$.  Therefore 
$g\in \alpha S$, which is a cotradiction.
Thus (\ref{eqn:2.1})  was verified.
By Lemma \ref{lemma:ord}, we have   
\begin{align*}
&~~~~~
 \alpha^{k} I^{*}(\omega, d\alpha) \in S_{(\alpha)}
\Leftrightarrow
 \mbox{\rm ord}_{\alpha} I^{*}(\omega, d\alpha) \le k\\
&\Leftrightarrow
\mbox{\rm ord}_{\alpha} I^{*}(\nabla_{D} (\omega), d\alpha)
= \mbox{\rm ord}_{\alpha} D(I^{*}(\omega, d\alpha)) \le k+2\\
&\Leftrightarrow
 \alpha^{k+2} I^{*}(\nabla_{D} (\omega), d\alpha) \in S_{(\alpha)},
\end{align*}
where $k := {\bf m} (\ker(\alpha))$.
This implies
\[
\omega\in\Omega(\A, {\bf m} )^{W}  \Leftrightarrow
\nabla_{D} (\omega)\in\Omega(\A, {\bf m}+2 )^{W}.
\]
 \owari

\section{Logarithmic derivation modules}
\label{sec:4} 
In this section, 
we translate our main results
Theorems \ref{thm:automorphism},
\ref{thm:primitivefiltration} and \ref{thm:bases} into the
corresponding theorems in the language of
logarithmic derivation modules.
Let $\Der_{S} $ and $\Der_{F} $ denote the 
$S$-module of $\R$-linear derivations 
from $S$ to itself and the $F$-vector space
of  $\R$-linear derivations 
from $F$ to itself.
Recall the $S$-linear isomorphism 
\[
I^{*} : \Omega_{F} \rightarrow \Der_{F},
\ \ \ \
I^{*}(\omega) (f) := I^{*}(\omega, df)  
\ \ \ \
(\omega\in \Omega_{F}, f\in F).
\]
The traslation of the main results is done
by
the isomorphism $I^{*} $.

\begin{define}
Define 
the $S$-module 
$D(\A,-\infty)$
of 
{\bf
logarithmic 
derivations}
by 
\begin{multline*}
D(\A,-\infty):
=\{\xi \in \Der_F ~|~Q^{N} \xi 
\mbox{~and~}
(Q/\alpha_{H})^{N} \xi(\beta)
\mbox{ ~are both regular}\\
\mbox{~for~}N \gg 0,
H 
\in \A 
\mbox{~and~}
\beta\in V^{*} 
\mbox{~with~}
 I^{*} (d\alpha_{H}, d\beta)=0
\}.
\end{multline*}
\end{define}
Then the map $I^{*} $ induces an $S$-linear isomorphism
\[
I^{*} :
\Omega(\A, \infty)
\widetilde{\longrightarrow}
D(\A, -\infty).
\]

Let $\A$ be a Coxeter arrangement which may not be irreducible
and $D$ be a primitive derivation.
The $T$-linear connection
\[
\nabla_{D} : \Der_{F} \rightarrow \Der_{F} 
\]
is characterized by
(1) 
$\nabla_{D} (f\xi) = D(f) \xi + f (\nabla_{D} \xi)$
and
(2) $\nabla_{D} (\partial_{v} ) = 0$ for all $v\in V$. 
Here the derivation $\partial_{v}  $ satisfies
$\partial_{v} (\alpha) = \alpha(v)$ for any $\alpha\in V^{*} $. 
Then it is not hard to see 
$
(\nabla_{D} \xi)(\alpha)
=
D(\xi(\alpha))
$ 
for all $\alpha\in V^{*} $.

\begin{lemma}
For $\omega\in \Omega_{F} $ we have
$$
I^{*} (\nabla_{D}(\omega))
=
\nabla_{D}(I^{*}(\omega)).
$$ 
In other words, the following diagram is commutative:

$$
\xymatrix@R1pc{
 \Omega(\A, \infty) \ar[r]^{\nabla_D} 
\ar[d]_{I^{*}} 
&\Omega(\A,  \infty) 
\ar[d]_{I^{*}} 
\\
D(\A, - \infty)    \ar[r]^{\nabla_D} 
&
D(\A, - \infty).
} 
$$
\end{lemma}

{\em Proof.}
It is enough to prove
$
I^{*} (\nabla_{D}(\omega))(\alpha)
=
\nabla_{D}(I^{*}(\omega))(\alpha)
$
for all $\alpha\in V^{*} $.
Compute
\[
(I^{*}(\nabla_{D} \omega) )(\alpha)
=
I^{*}(\nabla_{D} \omega, d\alpha)
=
D(I^{*}(\omega, d\alpha))
=
D(I^{*}(\omega)(\alpha))
=
(\nabla_{D} I^{*}(\omega))(\alpha).
\]
 \owari

\begin{define}
\label{def:primitivefiltrationD}
Let 
$\bf m : \A \rightarrow \Z$ be an arbitrary multiplicity.
Define
\begin{multline*} 
D(\A, {\bf m})
:=
I^{*} (\Omega(\A, -{\bf m} ))
=
\{
\xi\in D(\A, -\infty)
~|~
(Q/\alpha_{H})^{N} \xi(\alpha_{H}) \in \alpha_{H}^{{\bf m}(H) }  S \\
\mbox{~for all ~} H \in \A 
\mbox{~for~} N \gg 0 
\}
\end{multline*} 
and
$$
D(\A, {\bf m})^{W}
:=  
D(\A, {\bf m}) \cap 
D(\A, -\infty)^{W}.
$$
The {\bf
primitive filtration
of
$D(\A, -\infty)^{W}$
 induced from ${\bf m}$} 
is given by
\[
R^{({\bf m})}_{k} 
:=
\nabla_{D}^{k} 
D(\A, {\bf m})^{W}
  \ \ (k\in\Z).   
\]
\end{define} 
Note that 
$$
D(\A, {\bf m})
=
\{
\xi\in \Der_{S} 
~|~
\xi(\alpha_{H})
\in \alpha_{H}^{{\bf m}(H)}S
\mbox{~for all~} H\in\A\}  
$$
if ${\bf m}(H) \ge 0$ for all $H\in\A$.
In this case, $D(\A, {\bf m})$ was
introduced by Ziegler \cite{Z89}.  

\medskip

Theorem \ref{thm:primitivefiltration}
is translated into:

\begin{theorem}
The primitive filtration is an increasing filtration
\[
\dots \subset R^{({\bf m})}_{-1}  \subset R^{({\bf m})}_{0}
  \subset R^{({\bf m})}_{1} \subset \dots
\]
such that 
\[
R^{({\bf m})}_{k} = R^{({\bf m} - 2k)}_{0}  
=
D(\A, {\bf m} - 2k)^{W}. 
\]
\label{thm:primitivefiltrationD}
\end{theorem}

We construct bases for the primitive filtration of
$D(\A, -\infty)^{W}$
induced from 
$\bf 1$ 
by traslating Theorem \ref{thm:bases}
as follows:

\begin{theorem}
For any $k\in\Z$,
$1\le i\le t$  and $1\le j\le \ell[i]$ , define
$$
\xi_{j}^{(k)}[i]:=
\nabla_{D[i]}^{k} \left(I^{*} (dP_j[i])\right).
$$
Let
$$
\Xi^{(k)}[i] := \{\xi_{j}^{(k)}[i]\}_{1\le j\le \ell[i]}, 
\
\Xi^{(k)} := \bigcup_{i=1}^{t}  \Xi^{(k)}[i],
\ \mbox{~and~} \
\Xi := \bigcup_{k\in\Z} \Xi^{(k)}.
$$
Then
\begin{itemize}
\item[(1)]
the $S$-module
$
D(\A, -2k+1)
$
is free with a basis
$
\Xi^{(k)}, 
$ 
\item[(2)]
the $R$-module
$
D(\A, -2k+1)^{W}
$
is free with a basis
$
\Xi^{(k)}, 
$ 
\item[(3)]
the $T$-module
$D(\A,-2k+1)^{W} $
 is free with a basis
$
\bigcup_{p\le k}
\Xi^{(p)},  
$ and
\item[(4)]
the $T$-module
$D(\A, \infty)^{W} $
 is free with a basis
$
\Xi.  
$ 
\end{itemize}

\label{thm:basesD}
\end{theorem}

 \vspace{5mm}

\end{document}